\input amstex 
\documentstyle {amsppt} 
\vsize=7in \hsize=5.5in 
\topmatter \title The Weil-\'etale topology for number rings \endtitle 
\author
S.Lichtenbaum \endauthor

\thanks {  The author has been partially supported by the University of Paris and the Newton Institute} 
\endthanks 
 \endtopmatter
\document \magnification = \magstep 1 \baselineskip = 1.5\baselineskip

\overfullrule0pt
\define \Lim {\underrightarrow {Lim}}
\heading \S
0. Introduction \endheading
 
The purpose of this paper is to serve as the first step in the construction of a new Grothendieck topology (the Weil-\'etale
topology) for arithmetic schemes $X$ (schemes of finite type over Spec $\Bbb Z$), which should be in many ways better suited than
the
\'etale topology for the study of arithmetical invariants and of zeta-functions.  The Weil-\'etale cohomology groups 
of "motivic sheaves' or "motivic complexes of sheaves" shold be finitely generated abelian groups, and the special values of
zeta-functions should be very closely related to Euler characteristics of such cohomology groups. 

As an example of the above philosophy, let $\bar X$ be a compactification of $X$.  This involves first completing $X$ to obtain a
scheme $X_1$ such that $X$ is dense in $X_1$ and $f: X_1 \to$  Spec $\Bbb Z$ is proper over its image, and then, if $f$ is
dominant, adding  fibers over the missing points of Spec $\Bbb Z$ and the archimedean place of $\Bbb Q$ to obtain $\bar X$.

Let $\phi$ be the natural inclusion of $X$ into $\bar X$. The following should be true: 

a) The Weil-\'etale hypercohomology groups with compact support $H^q(\bar X, \phi_! \Bbb Z)$ are finitely generated abelian groups
which are equal to $0$ for all but finitely many $q$, and independent of the choice of compactification $\bar X$.  We will denote
them by
$H^q_c(X, \Bbb Z)$. 

b) If $\tilde \Bbb R$ denotes the "sheaf of real-valued functions" on $X$, then the cohomology groups $H^q(\bar X, \phi_!\tilde
\Bbb R)$ are independent of the compactification, and we denote them by $H^q_c(X, \tilde \Bbb R)$.  The natural map from
$H^q_c( X, \Bbb Z)\otimes \Bbb R$ to $H^q_c( X, \tilde \Bbb R)$ is an isomorphism.  (Note that this is not at all a
formality, and would for instance be false if we considered cohomology on all of $\bar X$).

c) There exists an element $\psi$ in $H^1(\bar X, \tilde \Bbb R)$ such that the complex $(H^*_c(X, \tilde \Bbb R), \cup
\psi)$ (this is a complex under Yoneda product with $\psi$) is exact. 

Then the Euler characteristic $\chi_c(X)$ of the complex $H^q_c(X, \Bbb Z)$ is well-defined (See Section 7), and we can
 describe the behavior of the zeta-function $\zeta_X(s)$ at $s = 0$ ($\zeta^*(0) = lim_{s \to 0} \zeta (s) s^{-a}$ where
$a$ is the order of the zero of $\zeta (s)$ at $s = 0$)  by the formula
$\zeta_X^*(0) = \pm \chi_c (X)$.   

Defining $\zeta^*(X, -n)$ in the analogous fashion, and taking advantage of the formula $\zeta_X(s) = \zeta_{X \times \Bbb A ^n} (s
+ n)$, we can conjecturally describe the behavior of the zeta-function of any arithmetic scheme at any non-positive integer $-n$ by
the formula
$\zeta^*(X, -n) = \pm \chi_c (X \times \Bbb A^n)$, 

There should exist motivic complexes $\Bbb Z(-n)$ whose Euler characteristics give the values of $\zeta^*(X, -n)$ directly, and
the above conjectural formula should give a guide to a possible definition.

In this paper we only define the Weil-\'etale topology in the case when $F$ is a global number field and $X =$ Spec $O_F$.  We then
compute the cohomology groups $H^q_c(X, \Bbb Z)$ for $q = 0. 1. 2, 3$, and verify that our conjectured formula holds true under
the assumption that the groups $H^q_c(X, \Bbb Z)$ are zero for $q > 3$.

It is not hard to guess possible extensions of the definition given here to arbitary $X$, once we have defined Weil groups and
Weil maps for higher-dimensional fields, both local and global.  Kato has made a very plausible suggestion of such a definition,
and we hope to return to this question in subsequent papers.

We close the introduction with two remarks: 

1) The definition given here would work for any open subscheme of a smooth projective
curve over a finite field.  Do the cohomology groups thus obtained agree with the ones defined in our earlier paper [L]?  This
seems highly likely, but we haven't checked it.

2) What is the relation of these conjectures to the celebrated Bloch-Kato conjectures? In general, they are not even about the
same objects.  The Bloch-Kato conjectures concern the Hasse-Weil zeta-function of a variety over a number field, and our
conjectures concern the scheme zeta-function of a scheme over Spec $\Bbb Z$.  If the scheme is smooth and proper over Spec $\Bbb
Z$, then the zeta-function of the scheme is the same as the Hasse-Weil zeta-function of the generic fiber, so then we can ask if
the conjectures are compatible.  Even this seems far from obvious, although presumably true.

\heading \S 1. Cohomology of topological groups \endheading

 Let $G$ be a topological group. We define a Grothendieck topology $T_G$ as follows: 

Let the category $Cat(T_G)$ be the category of $G$-spaces and $G$-morphisms. 
A collection of maps $\{\pi_i: X_i \to X\}$ will be called a covering (so an element of $Cov(T_G)$ if it admits local
sections: for every $x \in X$ there exists an open neighborhood $V$ of $x$, an index $i$ and a continuous map $s_i:V \to X_i$
such that
$\pi_i s_i = 1$. 

We verify easily that $Cat (T_G)$ has fibered products. It is immediate that $T_G$ satisfies the axioms for a
Grothendieck topology, and we call $T_G$ the "local-section topology". 

Let $A$ be a topological $G$-module. We define a presheaf of abelian groups $\tilde A$ on $T_G$ by putting $\tilde A(X) =
Map_G(X,A)$ (the set of continuous $G$-equivariant maps from $X$ to $A$).

\proclaim {Proposition 1.1} $\tilde A$ is a sheaf. \endproclaim

 Proof. We have to show $\tilde A$ verifies the sheaf axiom: Let $\{\pi_i:X_i \to X\}$ be a cover. Let $\theta_1$ and $\theta_2$
be the maps from $\prod Map_G(X_i, A)$ to
$\prod Map_G(X_i \times_X X_j, A)$ induced by the two projections, and let $\psi$ be the natural map from $Map_G (X,A)$ to $\prod
Map_G(X_i, A)$. We have to check that
 if $f$ is in $\prod Map_G(X_i, A)$ and $\theta_1(f) = \theta_2(f)$, there is a unique $g$ in $Map_G(X, A)$ such that
$f = \psi (g)$. 

It is clear that $g$ exists and is unique as a map of sets; we need only show that $g$ is continuous. This
follows immediately from the existence of local sections.

Define $C^p(G, A)$ to be $Map_G(G^{p+1}, A)$, where $G$ acts diagonally on $G^{p+1}$. Let $\delta_p$ map $C^p(G,A) $ to
$C^{p+1}(G, A)$ by the standard formula 

$$ \delta_p f(g_0, \dots g_{p+1}) = \sum^{p+1}_0 (-1)^i f(g_0, \dots \hat{g_i}, \dots g_{p+1}) $$ 

Then the cohomology $ H^p_c(G,A)$ of this complex is the continuous (homogeneous) cochain cohomology of $G$ with values in $A$.

Remark. By the usual computation, this cohomology is the same as the inhomogeneous continuous cohain complex of
$G$ with values in $A$.

 Let $*$ denote a point, with trivial $G$-action. 

\proclaim {Definition 1.2}. We define the cohomology groups $H^i(G, A)$ to be $H^i(T_G, *, \tilde A)$. \endproclaim 

\proclaim {Proposition 1.3} Let $0 \to A \to B \to C \to 0$ be an exact (as abelian groups) sequence of
$G$-maps of topological $G$-modules . Assume that the topology of $A$ is induced from that of $B$ and that the map from $B$ to
$C$ admits local sections. Then the sequence of sheaves on $T_G$ : $0 \to \tilde A \to \tilde B \to \tilde C \to 0$ is also
exact, and consequently there is a long exact sequence of cohomology $$ 0 \to H^0(G, A) \to H^0(G, B) \to H^0(G, C) \to H^1(G,A)
\dots$$ \endproclaim 

Proof. It is immediate that the sequence of sheaves is left exact. Let $X$ be a $G$-space and let $f: X \to
C$. Then the projection on the first factor makes the fibered product $X \times_C B$ a local section cover of $X$. Let $p_1$
and
$p_2$ be the projections from
$X\times_C B$ to $X$ and $B$ respectively, and $\lambda$ the map from $B$ to $C$. Then $p_1^*f = \lambda_* p_2$, so the map from
$\tilde B$ to $\tilde C$ is surjective. 

\proclaim {Proposition 1.4}. The \v Cech cohomology groups $\check H^p(*. \tilde A)$ =
$\check H^p(T_G, *, \tilde A)$ are functorially isomorphic to $H^p_c(G, A)$. \endproclaim 

Proof. By definition $\check H^p(*, \tilde A)$ is the direct limit of the groups $\check H^p(\Cal U, \tilde A)$, where $\Cal U$
runs through the set of coverings of $*$. It is immediate that the map from $G$ to $*$ is an initial object in the category of
covers, so $\check H^p(* , \tilde A) = \check H^p(\{G\}, \tilde A)$. But, by definition, this is the cohomology of the complex

$$Map_G(G, A) \to Map_G(G \times G, A) \to Map_G(G\times G \times G, A) \dots $$ 

which is just the definition of the homogeneous
continuous cochain complex. 

\proclaim {Corollary 1.5} Let $A$ be a $G$-module with trivial $G$-action. Then the cohomology
group $H^1(G, A)$ is naturally isomorphic to
$Hom_{cont}(G,A)$. \endproclaim 

Proof. In any Grothendieck topology, $H^1(F) = \check H^1(F)$ for any sheaf $F$. On the other
hand the continuous cochain cohomology group $H^1_c(G, A)$ ($G$ acting trivially) is well-known to be the group of continuous
homomorphisms from $G$ to $A$. \medskip 

Our next goal is to relate the cohomology of $G$ to the \v Cech cohomology of the
underlying topological spaces of $G$ and its products. 

\proclaim {Lemma 1.6} Let $G$ be a topological group and $X$ a topological
space. Let $W = G \times X$, and let $G$ act on $W$ by $g(h,w) = (gh, w)$. In the local section topology on $W$ every cover
$\{\pi_i:U_i \to W\}$ has a refinement by a cover of the form \{$G \times V_x$\}, where $V_x$ is a topological neighborhood of
the point $x$ in $X$. \endproclaim 

Proof. Let $x \in X$. There exists an open neighborhood $V_x$ of $x$, an open neighborhood
$T_x$ of the identity $e$ of $G$, an index $i$, and a section $\lambda_x: T_x \times V_x \to U_i$. Let $i_x$ be the inclusion of
$V_x$ in $X$. Define a map $\rho_x: G \times V_x \to U_i$ by $\rho_x(g,v) = g(\lambda_x(e, v))$. Clearly $\{G \times V_x\}$ is a
local section cover of $G \times X$.

 We have 

$$\pi_i\rho_x(g, v) = \pi_i g\lambda_x(e, v) = g\pi_i\lambda_x(e,v) = g(e, v) = (g, v)$$ 

This shows that $\pi_i\rho_x = id \times i_x$, and hence that $\{G \times V_x\}$ refines $\{U_i\}$. 

\proclaim {Corollary 1.7} a) Let $E$ be a local-section sheaf on $G \times X$. Define a local-section presheaf $\alpha_*E$ on $X$
by $\alpha_*E (Y) = E(G \times Y)$. Then $\alpha_*E$ is a sheaf for the local-section topology on $X$ and $\alpha_*$ is exact.

b) $H^q(T_G, G \times X, E)$ is isomorphic to $H^q(X, \alpha_* E)$.
\endproclaim 

We observe that we may restrict $\alpha_*F$ to the usual topology on $X$ and obtain the same cohomology, since
usual topological covers are cofinal in local-section covers. 

Proof.  a) Let $\{U_i \to Y\}$ be a local-section cover of $Y$. Then
$\{G \times U_i \to G \times Y\}$ is a local-section $G$-cover of $G \times Y$ and $G \times (U_i \times_Y U_j)$ is naturally
isomorphic to $(G \times U_i) \times_{G \times Y} (G \times U_j)$, so $\alpha_*F$ is a sheaf.

 We note that $\alpha_*$ is clearly
left exact. Let $E \to F$ be a surjective map of sheaves on $G \times X$. Let $x \in \alpha_*F(Y) =F(G \times Y)$. There exists
a local-section cover
$\{\pi_i:U_i \to G \times Y\}$ such that $\pi_i^*(x) \in F(U_i)$ lifts to $E(U_i)$. By Lemma 1.6, we may assume that $U_i$ is
$G \times V_i$, where the $V_i$'s are an open cover of $Y$, and $\pi_i = (id, \lambda_i)$, where $\lambda_i$ is the inclusion
of $V_i$ in $Y$. Clearly $\lambda_i^*x$ comes from $\alpha_*E(V_i)$, so $\alpha_*E \to \alpha_*F$ is surjective, and $\alpha_*$
is exact.

It is immediate that $\alpha^{-1}$, defined by $\alpha^{-1} (Y) = Y \times G$ is a map of topologies, so we have a Leray spectral
sequence for $\alpha_*$. This spectral seguence degenerates because $\alpha_*$ is exact, yielding the desired isomorphism.

In any Grothendieck topology, we have the presheaf $\Bbb Z'$, defined by assigning the group $\Bbb Z$ to any object amd the
identity to any map.  We define $\Bbb Z $ to be the sheaf associated with the presheaf $\Bbb Z'$.  If the topology
is
$T_G$, we also have the sheaf $\tilde \Bbb Z$, which corresponds to the trivial $G$-module $\Bbb Z$ and is characterized by
$\tilde \Bbb Z(X) = Map_G(X, \Bbb Z)$.  We can define a map from the presheaf $\Bbb Z$ to $\tilde \Bbb Z$ by sending $n$ to the
map with the constant value $n$.  This is clearly injective, and induces an injective map from $\Bbb Z$ to $\tilde \Bbb
Z$.

This map is also surjective. Let $f:X \to \Bbb Z$, and let $X_n = f^{-1}(n)$.  The $X_n$'s form a disjoint open cover of $X$, and
the Cech cohomology $H^0$ of the presheaf $\Bbb Z'$ with respect to this cover contains an element $g$ which is $n$ on $X_n$.  Then
$g$ determines an element of $\Bbb Z$ which maps onto $f$.

 \heading \S 2. An alternative definition \endheading

 Let $G$ be a topological group. We construct a simplicial $G$-space $S_n$ as follows: 

Let $S_n = G^{n+1}$, and let $G$ act on $S_n$ by $g(g_0, g_1 \dots g_n) = (gg_0, g_1 \dots, g_n)$. 

Now define face maps $\rho_i: S_n \to S_{n-1}$ by: $\rho_i(g_0,
\dots g_n) = (g_0, \dots, g_ig_{i+1}, \dots g_n)$ for $0 \leq i < n$, and $\rho_n(g_0, \dots g_n) = (g_0, \dots g_{n-1})$. 

The maps $\rho_i$ are maps of $G$-spaces, and a straightforward verification shows that $\rho_i\rho_j =\rho_{j-1}\rho_i$ if $i
< j$.

We will not use the degeneracy maps, so we omit the definition.

Now let $\tilde S_n = G^{n+1}$, but with $G$ acting diagonally on $\tilde S_n: g(g_0 \dots g_n) = (gg_0 \dots gg_n)$.  Let
$\pi_i:\tilde S_n \to \tilde S_{n-1}$ by $\pi_i(g_0 \dots g_n) = (g_0 \dots \hat g_i \dots g_n)$.   Computation shows a) $\pi_i$
is a $G$-map, and b) $\pi_i\pi_j = \pi_{j-1}\pi_i$ if $i < j$.

Let $\phi:\tilde S_n \to S_n$ by $\phi(g_0 \dots g_n)=(g_0, g_0^{-1}g_1, g_1^{-1}g_2 \dots g_{n-1}^{-1}g_n)$.

We verify that $\phi$ is a $G$-map and that $\rho_i\phi = \phi\pi_i$.

Let $F$ be a local-section sheaf on the site $T_G$.  Let $F_n$ be the sheaf on $G^n$ (as topological space) defined by $F_n(U)
= F(G \times U)$, where $G$ acts on $G \times U$ by acting by left translation on $G$ and trivially on $U$.

We define maps $\tilde \rho_i:G^n \to G^{n-1}$ by 

$$ \bar \rho_i(g_1 \dots g_n)= (g_1 \dots g_ig_{i+1} \dots g_n) ( 1 \le i \le n-1)$$

$$\bar \rho_0(g_1 \dots g_n) = (g_2 \dots g_n) $$

$$ \bar \rho_n(g_1 \dots g_n) = (g_1 \dots g_{n-1})$$

Let $p_n$ be the natural projection from $S_n = (G \times G^n)$ to $G^n$.  We check that $p_{n-1}\rho_i = \bar \rho_i p_n$, and
so automatically $\bar \rho_i \bar \rho_j = \bar \rho_{j-1}\bar \rho_i$.

Now take the (second) canonical flabby resolution $T_{j,n}$ of $F_n$ on $G^n$. [Some words are in order.  If $X$ is a topological
space and $F$ is a sheaf on $X$ the usual canonical flabby resolution is obtained by defining $C^0(F)$ to be $\prod_{x \in
X}(i_x)_*F_x$, embedding $F$ in $C^0(F)$, taking the quotient $G$,  embedding $G$ in $C^0(G)$, and continuing this process to
obtain a flabby resolution $0 \to F \to C^0(F) \to C^0(G) \to \dots $. On the other hand, the second canonical flabby resolution 
looks like (see [Go] \S 6.4 for details) $0 \to F \to C^0(F) \to C^0(C^0(F)) \to \dots $, after defining suitable coboundary
maps.  We have to use this construction to compare our definition with that of Wigner, who when he says "canonical flabby
resolution" means this one.] By construction we have for each
$i$ a map from
$F_{n-1}$ to $(\bar \rho_i)_* F_n$, which is easily seen to induce inductively a map from $T_{j,n}$ to $(\bar
\rho_i)_*T_{j,n}$, and hence a map from $\Gamma(G^{n-1}, T_{j,n-1})$ to $\Gamma (G^n, T_{j,n})$.  By taking the alternating
sum of these maps we get a map $\delta_{j,n} : \Gamma(G^{n-1},T_{j,n-1}) \to \Gamma(G^n, T_{j,n})$, and thus a double complex.
We define $\tilde H^*(G, F)$ to be the hypercohomology of this double complex.

\proclaim {Proposition 2.1}. The cohomology groups $H^i(T_G, *, F)$ are functorially isomorphic to the cohomology groups
$\tilde H^i(G,F)$. \endproclaim

Proof.  We need only check that $H^0 = \tilde H^0$, that the $\tilde H^i$'s form a cohomological functor, and that the $\tilde
H^i$'s vanish on injectives for $i > 0$.

 Since the canonical flabby resolution takes short exact sequences of sheaves into short exact sequences of complexes, the
$\tilde H^i$'s form a cohomological functor. (Recall that corollary 1.7 implies that an exact sequence of sheaves on $T_G$
gives rise to an exact sequence of sheaves on $G^n$ for every
$n$). 

If $F = I$ is injective, $I$ restricts to an injective sheaf $J_n$ on $G \times G_n$ (If $f$ is the map from $G \times
G_n$ to a point, $f^*$ takes injectives to injectives, since it has the exact left adjoint $f_!$). We know that $I_n =
\alpha_*J_n$ is flabby, hence acyclic, and so the homology of the flabby resolution of $I_n$ reduces to $H^0(G^n, I_n)$ and the
spectral sequence of a double complex shows that our hypercohomology is the cohomology of the complex $H^0(G_n, I_n) = I_n(G_n)
= I(G \times G_n)$. The equality $\rho_i\phi = \phi\pi_i$ shows that the homology of $I(G \times G_n)$ is the same as the
homology of $I(G^{n+1})$ with diagonal action, which is the \v Cech cohomology $\check H^i({G}, I)$. Since \v Cech cohomology
vanishes for injectives for $i > 0$ so does $\tilde H^i(G, F)$. 

Finally, it follows again from the formula $\rho_i\phi =
\phi\pi_i$ that if $F$ is any presheaf on the category of $G$ -spaces. the cohomology of the complex $F(S_n)$ is naturally
isomorphic to the cohomology of the complex $\alpha_*F(G^n) = F(\tilde S_n)$, where the coboundary maps are the alternating
sums of the maps induced by $\rho_i$ and $\pi_i$ respectively. It follows that the cohomology $\tilde H^0(G, F)$ is naturally
isomorphic to the \v Cech cohomology $\check H^0(G, F)$ which in turn is $H^0(T_G, *, F)$. 

{\bf Remark 2.2} We observe that if
$F$ is a sheaf of the form $\tilde A$, then our cohomology groups are exactly the cohomology groups denoted by $\hat H^*(G, A)$
by David Wigner ([W], p.91). We then obtain as a corollary of Wigner's Theorem 2 ([W], p.91) that if $G$ is locally
compact, $\sigma$-compact, finite dimensional, and $A$ is separable and has Wigner's "property F", then our $H^*(G, A)$ are
naturally isomorphic to the groups $H^*(G, A)$ defined by Wigner in [W], (which we will call $H^*_{Wig}(G,A)$). We further
point out that under the same conditions Wigner's groups are naturally isomorphic to the groups (which we wil call $H^*_M(G,
A)$) defined by Calvin Moore in [M] and used by C.S. Rajan in [R]. (Wigner's Theorem 2 does not explicitly require
separability, but his proof that certain categories of modules are quasi-abelian is not valid without it.) In order to apply
this result, we recall that Proposition 3 of [W] tells us that any locally connected complete metric topological group
(for instance, $\Bbb Z$, $S^1$, or $\Bbb R$) has property F.

 \proclaim{ Theorem 2.3} There is a spectral sequence $E^{p,q}_1 =
H^q_{top}(G^{p}, \alpha_*F) \Rightarrow H^{p+q}(T_G,*, F)$. \endproclaim 

Proof. This is just the spectral sequence of the
double complex defining $\tilde H^*(G,F)$. 

\proclaim {Corollary 2.4} Let $G$ be a) a profinite group, or b) the Weil group of
a global function field , and $A$ a topological $G$-module. Then the cohomology groups $H^i(G, A)$ are canonically isomorphic
to the usual groups $H^i_{cont}(G, A)$ given by complexes of continuous cochains. \endproclaim 

Proof. We show first that the
cohomological dimension of a profinite space $X$ is zero. To do this it suffices to show (by using alternating cochains)) that
every open cover has a refinement by a disjoint cover. It is immediate that $X$ has a base for its topology consisting of sets
$U_i$ which are both open and closed. By compactness, any cover has a refinement $\{U_1, \dots U_n \}$ consisting of finitely
many such $U_i$ . Let $C(U) = X-U$. Then $ \{U_1, U_2 \cap C(U_1), U_3\cap C(U_2) \cap C(U_3) \dots \}$ is a further refinement
which is disjoint. 

In case a) each $G^q$ is profinite, so has cohomological dimension zero, and in case b) the Weil group $G$
is the topological product of a profinite group and a discrete group, so $G^q$ is the disjoint union of open profinite spaces,
so again has cohomological dimension zero. So in each case the spectral sequence degenerates, to yield that $H^*(G, F)$ is the
cohomology of the complex $F(G \times G^p)$. We have seen that this is the same as the cohomology of the complex $F(G^{p+1})$,
with $G$ acting diagonally, which is just the homogeneous continuous cochain complex of the $G$-module $A$ if $F= \tilde A$.

\proclaim {Lemma 2.5} Let $X$ be the product of a compact space and a metrizable space, and let $E$ be a sheaf of modules over 
the sheaf of continuous real-valued functions on $X$. Then $H^q(X^p, E) = 0$ for all $p, q > 0$.
\endproclaim 

Proof. The hypothesis implies
that for any $p$, $X^p$ is again the product of a compact space and a metrizable space, and so paracompact. We recall from
([Go], p.157) that any sheaf of modules over the sheaf of continuous real-valued functions on a paracompact space is fine, so
"mou", so acyclic.

\proclaim {Corollary 2.6} Let $G$ be a topological group which is, as topological space, the product of a compact space and a
metrizable space (e.g. the Weil group of a global or local field) and let $\Bbb R$ denote the real numbers with trivial
$G$-action. Then the cohomology groups $H^p(G, \tilde \Bbb R)$ are given by the cohomology of the complex of homogeneous
continuous cochains from $G$ to $\Bbb R$. \endproclaim 

Proof. Let $F = \Bbb R$. We first observe that $\alpha_*(F) (U) = F(G
\times U) = Map_G(G \times U, \Bbb R)$, which is naturally isomorphic to $Map (U, \Bbb R) = \tilde \Bbb R(U)$. Then Lemma 2.5
implies that the spectral sequence of Theorem 2.3 degenerates, so that the cohomology $H^p(G,\tilde \Bbb R)$ is given by the
cohomology of the complex $H^0(G^p, \alpha_* \Bbb R) = Map_G(G \times G^p, \Bbb R) = Map_G(S_p,\Bbb R)$. As above, this is the
same as the cohomology of the homogeneous cochain complex $Map_G(\tilde S_p, \Bbb R)$. 

\heading \S 3. Cohomology of the Weil group. \endheading

Let $F$ be a number field (resp. a local field), $\bar F$ an algebraic closure of $F$, and $G_F$ the galois group of $\bar F$ over
$F$. Let
$K$ be a finite Galois extension of
$F$. and let $C_K$ denote the id\`ele class group of $K$ (resp. $K^*$).

 Now fix a Weil group $W_F$ associated with the topological class formation $Lim (C_K)$, where the limit is taken
over fields $K$ finite and Galois over $F$. We recall that $W_F$ is equipped with a continuous homomorphism $g: W_F \to G_F$.
If $K$ is such a field, let $W_K = g^{-1}(G_K)$, and let $W_K^c$ be the closure of the commutator subgroup
of $W_K$ in $W_F$. Then it is shown in [Artin-Tate ] that $W_F/W_K^c$ is a Weil group for the pair $(G(K/F), C_K)$ (resp.
$(G(K/F), K^*)$). So having fixed a Weil group $W_F$, we have canonical maps from it to $W_{K/F} = W_F/W_K^c$. The standard
construction of the Weil group
$W_F$ (See [A-T]) shows that $W_F$ is the projective limit of the groups $W_{K/F}$.

Now let $F$ be global and $S$ be a finite set of valuations of $F$ including the archimedean valuations and all valuations which
ramify in
$K$, but not including the trivial valuation.  Let $U_{K,S}$ be the subgroup of the id\`ele group $I_K$ consisting of those
id\`eles which are 1 at valuations lying over  $S$, and units at valuations not lying over $S$.  It is well known (see [N]. p.
393) that
$U_{K,S}$ is a cohomologically trivial $G(K/F)$-module. The natural map from $U_{K,S}$ to the id\`ele class group $C_K$ is
obviously injective and we identify $U_{K,S}$ with its image.  Let the $S$-id\`ele class group $C_{K,S}$ be defined by $C_{K,S} =
C_K/U_{K,S}$.  Then the natural maps from the Tate cohomology groups $\hat H^i(G(K/F), C_K)$ to $\hat H^i(G(K/F),C_{K,S})$ are
isomorphisms for all $i$.

Let $\alpha$ be the fundamental classs in $\hat H^2(G(K/F), C_K)$ and $\beta$ its image in $\hat H^2(G(K/F), C_{K,S})$.  It follows
immediately from the fact that $C_K$ is a class formation that. for all $i$. cup-product with $\beta$ induces an isomorphism
between
$\hat H^i(G(K/F),\Bbb Z)$ and $\hat H^{i+2}(G(K/F), C_{K,S})$.  We then define the $S$-Weil group $W_{K/F,S}$ to  be the extension
of
$G(K/F)$ by $C_{K,S}$ determined by $\beta$,  There is clearly a natural surjection $p_S$ from $W_{K/F}$ to $W_{K/F,S}$, and it
follows from the arguments in [A-T] (p.238) that there is a natural isomorphism from $W_{K/F,S}^{ab}$ to $C_{F,S}$.

  Let $N_{K,S}$ be the kernel of the natural map from $W_F$ to $W_{K/F,S}$.  Let $A$ be a topological $W_F$-module, and let
$A_{K,S}$ be the topological
$W_{K/F,S}$ module consisting of the invariant elements $A^{N_{K,S}} \subseteq A$.  Assume that $A = \bigcup A_{K,S}$.

\proclaim{Lemma 3.1} The Weil group $W_F$ is the projective limit over $K$ and $S$ of the groups $W_{K/F, S}$. \endproclaim

Proof.  It suffices to show that the relative Weil group $W_{K/F}$ is the projective limit over $S$ of the groups $W_{K/F,S}$. 
The maps $p_S$ induce a map $p$ from $W_{K/F}$ to the projective limit.  Let $W^1_{K/F}$ (resp.$ W^1_{K/F, S})$ be the kernel of
the absolute value map of $W_{K/F}$ (resp.$W_{K/F, S})$ to $\Bbb R^*$.  Since $W^1_{K/F}$ is compact and the maps $p_S$ are
surjective, $p$ is surjective as a map from $W^1_{K/F}$ to the projective limit of $W^1_{K/F, S}$, and hence $p$ is surjective. 
The proof that $p$ is injective immediately reduces to showing that the map from $C_K$ to the projective limit of the
$C_{K/F,S}$ is injective, which in turn follows from the corresponding fact for the id\`ele groups. 
 
 \proclaim {Definition 3.2}. We define the cohomology group $H^q(W_F, A)$ to be the direct limit of the cohomology groups
$H^q(W_{K/F,S}, A_{K,S})$. \endproclaim

The cohomology groups of $W_{K/F,S}$-modules are the ones defined in Section 1. We observe that $W_{K/F,S}$ is locally compact,
$\sigma$-compact and finite-dimensional, so Wigner's comparison theorem applies and the cohomology of these groups with
coefficients in $\Bbb Z$, $\Bbb R$ or $S^1$ are the same as Wigner's cohomology groups and therefore also Moore's cohomology
groups.

We now compute the cohomology groups $H^q(W_F, \Bbb Z)$: Evidently $H^0(W_F, \Bbb Z) =\Bbb Z$. Since $H^1(W_{K/F,S},
\Bbb Z)= Hom_{cont}(W_{K/F,S}, \Bbb Z) = 0$, (because $W_{K/F,S}$ is an extension of a compact group by a connected group), we
have $H^1(W_F, \Bbb Z) = 0$.

 We have the following result of Moore
([M2], Theorem 9, p.29), as quoted by Rajan ([R], Proposition 5):

 \proclaim {Lemma 3.3} Let $G$ be a locally compact
group. Let $N$ be a closed normal subgroup of $G$ and let $A$ be a locally compact, complete metrizable topological
$G$-module. Then there is a spectral sequence $$ E^{p,q}_2 \Rightarrow H^{p+q}_M(G,A) $$ where $E^{p.q}_2 = H^p_M(G/N,
H^q_M(N,A))$ if $q = 0, q = 1$, or $ H^q_M(N, A) = 0$ \endproclaim 

\proclaim {Lemma 3.4} The cohomology groups $H^q(W_{K/F,S}, \Bbb R)$ are: $\Bbb R$ if $q = 0$,
$Hom_{cont} (\Bbb R, \Bbb R) \simeq \Bbb R$ if $q=1$, and $0$ if $q > 1$. \endproclaim

 Proof. We know by Corollary 2.5 that
these cohomology groups are given by the continuous cochain cohomology. It is well-known (See [B-W]) that if $G$ is compact then
$H^q(G, \Bbb R) = 0$ for $q > 0$, and $H^0(\Bbb R, \Bbb R) = \Bbb R$, $H^1(\Bbb R, \Bbb R) = Hom_{cont}(\Bbb R, \Bbb R)$ , and
$H^q(\Bbb R. \Bbb R) = 0$ for $q >1$. The result then follows from the fact that we have the exact sequence: 

$$ 1 \to W^1_{K/F,S} \to W_{K/F,S} \to \Bbb R \to 1 $$

 with $W^1_{K/F,S}$ compact. and applying Lemma 3.3.

\proclaim {Lemma 3.5} $H^q_M(\Bbb R, \Bbb Z) = 0$ for $q > 0$.
\endproclaim Proof. $H^q_M(\Bbb R,
\Bbb Z) = H^q_{Wig}(\Bbb R, \Bbb Z) =$ (by Theorem 4 of [W]) $H^q(B_{\Bbb R}, \Bbb Z) = 0$ because $\Bbb R$ is
contractible. 

We see from Lemma 3.4 and the exact sequence $0\to \Bbb Z \to \Bbb R \to S^1 \to 0$ that we have $0 \to
Hom_{cont}(\Bbb R, \Bbb R) \to H^1(W_{K/F,S},S^1) \to H^2(W_{K/F,S}, \Bbb Z) \to 0$. Since the abelianized Weil group
$(W_{K/F,S})^{ab}$ is naturally isomorphic to $C_{F,S}$, $H^1(W_{K/F,S},S^1)$ is the Pontriagin dual $C_{F,S}^D$ of $C_{F,S}$,
which yields that $H^2(W_{K/F,S}, \Bbb Z)$ is the Pontriagin dual $(C^1_{F,S})^D$ of the id\' ele class group of norm one.  By
taking limits over $K$ and $S$ we obtain that $H^2(W_F, \Bbb Z) = (C_F^1)^D = 0$.

We next wish to show that $H^3(W_F, \Bbb Z) = 0$, and to do this it is, by Lemma 3.4, enough to show that $H^2(W_F, S^1) = 0$.
We first observe that Rajan's proof in [R] that the Moore cohomology groups $H^2_M(W_F, S^1)= 0$ works equally well to show
that $H^2_M(W^1_F, S^1) = 0$. Since for Moore cohomology, the cohomology of the projective limit of compact groups is the
direct limit of the cohomologies, ([M] or [R]) we have that $0 = H^2_M(W^1_F, S^1) = \Lim H^2_M(W^1_{K/F,S}, S^1) = $ (by
Remark 2.2) $\Lim H^2(W^1_{K/F,S},S^1) $. whixh is in turn equal to (by Lemma 3.4) $\Lim H^3(W^1_{K/F,S}, \Bbb Z)$.

It is easy to see that the Weil group
$W_{K/F,S}$ is the direct product (in both the algebraic and topological senses) of $W^1_{K/F,S}$ and $\Bbb R$. Applying the
Hochschild-Serre spectral sequence (Lemma 3.3) coming from the exact sequence $1 \to \Bbb R \to W_{K/F,S} \to W^1_{K/F,S} \to
1$, and using Lemma 3.5, we conclude that $H^q(W_{K/F,S}, \Bbb Z) = H^q(W^1_{K/F,S}, \Bbb Z)$. So $H^3(W_F, \Bbb Z)$ = (by
definition) $\Lim H^3(W_{K/F,S}, \Bbb Z) =\Lim H^3(W^1_{K/F,S}, \Bbb Z) = 0$. We sum up what we have shown in the following
theorem: 

\proclaim {Theorem 3.6} The cohomology groups $H^q(W_F, \Bbb Z)$ are given by: $H^0(W_F, \Bbb Z)= \Bbb Z$,\break
$H^1(W_F, \Bbb Z) = 0$, $H^2(W_F, \Bbb Z) = (C^1_F)^D $(the Pontriagin dual of $ C^1_F$) and $ H^3(W_F, \Bbb Z) = 0$.
\endproclaim 

Unfortunately so far we have not succeeded in computing the cohomology groups $H^q(W_F, \Bbb Z)$ for $q > 3$.

\proclaim{Lemma 3.7} If $v$ is not in $S$, the natural map induced by $\theta_v$ from $W_v$ to $W_{K/F,S}$ annihilates the
kernel $I_v$ of the natural map from $W_{F_v}$ to $\Bbb Z$.
\endproclaim

Proof.  Let $w$ be the valuation lying over $v$ determined by $\theta_v$.  The following diagram is commutative:

$$ \CD 1 @ >>> K_w^* @ >>> W_{K_w/F_v} @ >>> G(K_w/F_v) @>>> 1  \\
@ .@ VwVV @ V\pi_v VV @ViVV  \\
1@>>> \Bbb Z @ >f>> \Bbb Z @ >>> G(\kappa (w)/\kappa (v)) @ >>> 1 \\
\endCD $$
because fundamental classes of the two extensions correspond. (Here $i$ is the natural isomorphism and $f$ is the residue field
degree).  It follows that the image of $I_v$ in $W_{K_w/F_v}$ is isomorphic to the unit group $Ker(w$), which goes to zero in
$C_{K/F,S}$ and so a fortiori in
$W_{K/F,S}$.

\heading \S 4. The global Weil-\' etale topology \endheading 

Let $F$ be a global field     choose an algebraic closure $\bar F$ of $F$.
Let $G_F =G(\bar F/F)$ be the Galois group of $\bar F/F$.

 Let $v$ be a valuation of $F$, and $F_v$ the completion of $F$ at $v$. Choose an algebraic
closure $\bar F_v$ of $F_v$, and an embedding of $\bar F$ in $\bar F_v$. Choose a global Weil group $W_F$
and a local Weil group $W_{F_v}$. For each finite extension $E$ of $F$ in $\bar F$, let $E_v = EF_v$ be the induced completion
of $E$. Let $w$ be a valuation of $\bar F$ lying over $v$, and let $i_w^*$ be the natural inclusion of $G_{F_v}$ in $G_F$
whose image is the decomposition group of $w$. 

\proclaim {Definition 4.1} A Weil map $\theta_v$ is a continuous homomorphism from
$W_{F_v}$ to $W_F$ such that there exists a valuation $w$ of $\bar F$ such that the following diagrams are commutative for all
finite extension fields
$E$ of $F$: 

$$ \CD W_{F_v} @ >>> G_{F_v} @. \qquad \qquad E_v^* @ >>> W^{ab}_{E_v} \\ 
@ V\theta_v VV @ VV
i_w^*V \qquad \qquad @V n_v VV @ VVV \\ 
W_F @ >>> G_F @.\qquad \qquad C_E @ >>> W_E^{ab} \\ \endCD $$ 
where $n_v$
maps $a \in E^*$ to the class of the id\' ele whose $v$-component is $a$ and whose other components are $1$, and the map from
$W_{E_v}^{ab}$ to $W_{E}^{ab}$ is induced by $\theta_v$. \endproclaim 

It is an easy consequence of [T] that Weil maps
always exist, and are unique up to an inner automorphism of $W_F$.  

The local Weil group $W_v=
W_{F_v}$ maps to ${W_v}^{ab}$ = $F_v^*$, which in turn maps to $\Bbb Z$ by the valuation map $v$. Let $I_v$ be the kernel of
the composite map from $W_v$ to $\Bbb Z$. 

We choose once and for all a set of Weil maps $\theta_v: W_v \to W_{\Bbb Q}$ for all valuations $v$ of $\Bbb Q$.  If
$w$ is any valuation of a number field $F$, the inclusion of $W_w$ in $W_v$ and $\theta_v$ induce a Weil map $\theta_w:W_w \to
W_F$.

 Let $\bar Y = \bar Y_F$ be the set of all valuations of $F$.  We require the trivial valuation $v_0$ to be in $\bar Y$,
corresponding to the generic point of Spec $O_F$, where $O_F$ is the ring of integers of $F$>.  Let $W_{\kappa(v)}$ be $\Bbb
Z$ if $v$ is non-archimedean, $\Bbb R$ if $v$ is archimedean. and $W_F$ if $v = v_0$. We say that $v$ is a specialization of $w$ if
$w$ is $v_0$ and $v$ is not. In each case there is a natural map $\pi_v$ from the local Weil group $W_v$ to $W_{\kappa (v)}$, and
we let $I_v$ be its kernel.  It is an easy exercise to verify that if $K_w$ is a finite Galois extension of $F_v$, then the
map $\pi_v$ factors through $W_{K_w/F_v}$.

Let $K$ be a finite Galois extension of $F$. Let $S$ be a finite set of non-trivial valuations of $F$, containing all the
valuations of $F$ which ramify in $K$.   We now define a Grothendieck topology
$T_{K,S,\bar Y}$:

We first define a category Cat $T_{K,S,\bar Y}$ The objects of Cat $T_{K,S, \bar Y}$ are collections $((X_v),(f_v))$, where
 $v$ runs through all points of $\bar Y$, $X_v$ is a $W_{\kappa(v)}$-space,  and
if $v$ is a specialization of $w$, $f_v:X_v \to X_{w}$ is a map of $W_v$-spaces. (We regard $X_v$ as a $W_v$-space via
$\pi_v$, and $X_w$ as a
$W_v$-space via the Weil map $\theta_v$).  If $v = v_0$, we require that the action of $W_F$ on $X_v$ factor through
$W_{K/F, S}$.    

A morphism $g$ from $\Cal X =((X_v),(f_v))$ to $\Cal X' = ((X'_v),(f'_v))$ is a collection of $W_v$-maps $g_v: X_v \to
X'_v$ such that  $g_{v_0}f_v = f'_vg_v$.

We say that $g$ is a local section morphism if the maps $g_v$ from $X_v$ to $g_v(X_v)$ admit local sections.

The fibered product of $((X_{1,v}), (g_{1,v}))$ and $((X_{2,v}), (g_{2, v}))$ over 
$((X_{3,v}), (g_{3,v}))$ is given by  \linebreak $((X_{1,v} \times_{X_{3,v}}
X_{2, v}),((g_{1,v} \times g_{2,v}))$.

We define the coverings Cov $(T_{K,S, \bar Y})$ by:

A family of morphisms in our category $\{((X_{i,v}),(f_{i,v})) \to ((X_v),(f_v))\}$ is a
cover if $\{X_{i,v} \to X_v\}$ is a local section cover for all
$v$. 

Our category clearly has a final object $*_{(K,S)}$ whose components are the one-point space for each $v$ in $\bar Y$.

If $E$ is a sheaf for our topology, we define $H^i(\bar Y_{K,S}, E)$ to be $H^i(T_{K,S, \bar Y}, *_{(K,S)}, E)$, 
 and $H^i(\bar Y, E)$ to be the direct limit over $K$ and $S$ of the $H^i(\bar Y_{K,S},
E)$.

  We define a morphism of topologies $i_v^{-1}$ from $T_{K,S}$ to $T_{W_{\kappa(v)}}$ by $i_v^{-1}((X_v), (f_v)) = X_v$.  We have
the corresponding direct image maps $(i_v)_*$ from sheaves on $T_{W_{\kappa(v)}}$ to sheaves on $T_{K,S}$ by $(i_w)_*(E)
((X_v), (f_v)) = E(X_w)$.  For psychological reasons we define $j$ to be $i_{v_0}$.   It is clear that $i_v^{-1}$ preserves covers
and fibered products, and so is a morphism of topologies.

\proclaim {Definition 4.2} Let $\theta:H \to G$ be a morphism of topological groups, and $X$ an $H$-space.  Define $X \times^H
G$ to be the quotient (with the quotient topology) of $X \times G$ by the equivalence relation $(x, g) \sim (x', g')$ iff there
exists a $\tau \in H$ such that $x' = \tau x$ and $g' = g\tau^{-1}$.\endproclaim

Remark.  The functor which takes an $H$-space $X$ to the $G$-space $X \times ^H G$ is easily seen to be left adjoint to the
forgetful functor  from $G$-spaces to $H$-spaces, regarding a $G$-space as an $H$-space via $\theta$.

\proclaim {Lemma 4.3} Let $G$ be a topological group and let $I$ be a closed subgroup such that the projection $\rho$ from
$G$ to $G/I$ admits local sections.  Then the  category of
$G$-spaces with maps to $G/I$ is equivalent to the category of $I$-spaces and the covers in the respective categories
correspond.
\endproclaim

Proof.  If $X$ is an $I$-space let $\alpha(X) = (X \times ^I G, \lambda)$, 
where $\lambda: X \times ^I G \to G/I$ is
given by $\lambda(x, \sigma) =$ the coset $\sigma I$. 

 If $Z$ is a $G$-space with a map $\pi: Z \to G/I$, let $\beta(Z, \pi) = \pi^{-1}(I)$.  It is straightforward to verify that
$\alpha$ and $\beta$ are inverse functors. 


We now claim that the covers correspond. 

\proclaim {Lemma 4.4} If $\rho:G \to G/I$ admits local sections and the cover $\{X_i \to X\}$ admits local sections, then the
cover
$\{X_i \times ^I G \to X \times ^I G\}$ admits local sections. \endproclaim 

Proof.  Let $y = [x, \sigma]$ be the class of $(x, \sigma)$ in $X \times ^I G$.  Let $U$ be a neighborhood of $\rho(\sigma)$
such that there exists a continuous section $s: U  \to G$ of $\rho$.  Let $V = \rho^{-1}(U)$. Let $U^* = s(U)$.   We claim that
$X \times ^I V $ is functorially isomorphic to $X \times U^*$.  It is immediate that given $[x, v]$ in  $X \times ^I V$, there
exists a unique pair $(x', v')$ in $X \times U^*$ such that $[x, v] = [x', v']$.  In fact $(x', v') = ((s\rho(v))^{-1}vx,
s\rho(v))$.

So if $\{X_i \to X\}$ admits local sections so does $\{X_i \times U^* \to X \times U^*\}$, and then so does $\{X_i \times ^I V
\to X \times ^I V\}$, and therefore also $\{X_i \times ^I G \to X \times ^I G\}$.


\proclaim {Lemma 4.5}. If $I$ is a locally compact subgroup of a Hausdorff topological group $G$, the natural projection from
$G$ to $G/I$ is a fibration, and hence admits local sections. \endproclaim

Proof.  This is [W], Proposition 2, p.88. 

So we have proved 

\proclaim {Theorem 4.6}.  Let $G$ be a Hausdorff topological group, $I$ a locally compact subgroup, and $A$ a continuous
$G$-module.  Then $H^i(T_G, G/I,  \tilde A)$ is naturally isomorphic to $H^i(I, A)$.
\endproclaim

\proclaim {Theorem 4.7} Let $j = j_{\bar Y}$, and let $A$ be a topological $W_F$-module.  There exists a spectral sequence $$
E^{p.q}_2 = H^p(\bar Y, R^qj_*\tilde A)
\Rightarrow H^{p+q}(W_F, A)$$
\endproclaim 

Proof.  This follows from [A]( p. 44) by applying his Theorem 4.11 to $j = j_{K,S}$ and taking direct limits 
over $K$ and $S$.

The rest of this section will be devoted to computing the sheaves $R^qj_* \tilde A$. 
 Let $v$ be in $\bar Y$  Our goal is to prove:

\proclaim {Theorem 4.8}. Let $q > 0$, and let $B =B_q$ be $R^q(j_{K,S})_* \tilde A$. Then the natural map
from
$B$ to
$\coprod_{v\in S}{i_v}_*i_v^*B$ given by adjointness is an isomorphism of sheaves.
\endproclaim

We begin with:

\proclaim {Lemma 4.9}. Let $E$ be a Weil-\'etale sheaf on $T_{K.S.\bar Y}$.  Then $i_v^*E = 0$ for all $v \in \bar Y$ implies that
$E = 0$. \endproclaim

Proof.  We know that $i_v^*E$ is the sheafification of the presheaf inverse image $i_v^pE$.  If $X_v$ is a
$W_{\kappa(v)}$-space,
$i_v^pE(X)$ is the direct limit of $E(U)$, where $U = ((X'_v),(f_v))$ is an object of Cat $(T_K)$ such that there is a
map from $X_v$ to $i_v^{-1}(U) = X'_v$.  Since there exists a $U$ (for example the object which has $X_v$ at $v$, $X_v \times
^{W_v} W_{K/F,S}$ at the generic point, and the empty set elsewhere) with $i_v^{-1}(U) = X_v$ we may always assume that
$i_v^{-1}(U) = X_v$, i.e. that $X'_v = X_v$.  

More generally, if $h_v: Z_v \to X_v$ is a map of $W_{\kappa(v)}$-spaces, and $f_v: X_v \to X_{v_0}$ is a map of $W_v$-spaces,
then the map $h'_v: X'_{v_0} = Z_v \times ^{W_v} W_{K/F,S} \to X_{v_0}$ given by $h'_v(z, w) = wf_vh_v(z)$ is well-defined and so
we get a map of $(X'_{v_0}, Z_v, \phi,... \phi)$ to $(X_{v_0}, X_v, (X_w))$ which induces the original map $h_v$.  

If $U = (X_v)$ ia an object of Cat $T_{K,S, \bar Y}$, and $\alpha \in E(U)$, there is a covering $\{X_{v_i}\}$ of $X_v$ such that
$\alpha$ goes to zero in  each $i^pE(X_{v_i})$ = $E(U_i)$, where the $v$-component of $U_i = X_{v,i}$.  By the argument in the
preceding paragraph, we can induce these coverings from families of maps to $U$, and the collection of all these families will
be a covering of $U$ in which $\alpha$ goes to zero, thus making $\alpha = 0$.

\proclaim {Lemma 4.10} a)$i_v^*$ is exact. b) $i_v^*{i_v}_*i_v^*E$ is canonically isomorphic to $i_v^*E$. c)
$i_w^*{i_v}_* = 0$ if $v \neq w$. d) ${i_v}_*$ is exact.\endproclaim

Proof.  a) Since $i_v^*$ is a left adjoint, it is right exact.  Suppose that the sheaf $E$ injects into the sheaf $E'$, and
that $\alpha \in i_v^*E(X_v)$ goes to zero in $i_v^*E'(X_v)$. There exists a Weil-\'etale cover $(X_{v,i})$ of
$X_v$ and objects
$\Cal X_i$ of Cat $T_{K,S, \bar Y}$ such that for each $i$, $\alpha$ restricted to $X_{v,i}$ comes from  an element $\beta_i$ in
$E(\Cal X_i)$,
$(\Cal X_i)_v = X_{v,i}$, and the image of $\beta_i$ in $E'(\Cal X_i)$ is equal to zero.  Hence $\beta_i = 0$ and since $\alpha$
goes to zero in a Weil-\'etale cover, we have
$\alpha = 0$.   

b) This is a formal consequence of the fact that $i_v^*$ is left adjoint to ${i_v}_*$.

c) $(i_w^pF)(X_w)$ is the direct limit of $F(U =(X_{w_0}, X_w, \phi, .....\phi))$ where $X_w \to X_{w_0}$. If $F = (i_v)_*E$ and
$v \neq w$ then $F(U) = E(\phi) = 0$.

 d) This follows immediately from the fact that $v$ is a specialization of $v_0$, and if
$\Cal X = (X_{v_0}, X_v, (X_w (w
\neq v, v_0))$, then any covering
$X_{v,i}$ of
$X_v = i_v^{-1} (\Cal X)$ comes from the covering $\Cal X_i = (X_{v_0}, X_{v_i}, (X_w))$ of $\Cal X$.

\proclaim {Lemma 4.11} Let $E$ now be the sheaf $R^q(j_{K,S})_* \tilde A$, with $q > 0$.  If $v$ is not in $S$, then $i_v^*E = 0$.
\endproclaim

Proof.  Given a $W_{\kappa(v)}$ - space $X_v$ and an element $\alpha$ in $i_v^p(E)( X_v)$, we wil produce a cover $\{X_{v,i}\}$ of
$X_v$ such that the restriction of $\alpha$  vanishes on each $X_{v,i}$.  By Lemma 3.7, if $v$ is not in $S$, the Weil map
$\theta_v$ from $W_v$ to $W_{K/F,S}$ factors through $W_{\kappa(v)}$. So let us define  $X_{v_0}$ to be $ X_v \times
^{W_{\kappa(v)}} W_{K/F,S}$.   By using the definition of $i_v^p$, , $i_v^p(E)(X_v)$  is easily seen to be $E(X_{v_0}, X_v, \dots,
\phi, \dots )$, where all the spaces $X_w$ for $ w \neq v, v_0$ are empty.  By passing to a cover, we may assume that $\alpha$
comes from an element $\beta$ is $H^q(X_{v_0}, \tilde A)$. Since $q> 0$ and higher cohomology dies in a cover, we may choose a
cover
$X_{v_0, i}$ of $X_{v_0}$ such that $\beta$ goes to zero in $H^q(X_{v_0,i}, \tilde A)$. Letting $X_{v, i} = X_{v_0}
\times_{X_{v_0}} X_v$, we see that $\alpha$ goes to zero on each $X_{v,i}$.

Proof of Theorem 4.8: By lemma 4.11, $\prod_{v \in \bar Y }(i_v)_*i_v^*B$ is equal to $\prod_{v \in S} (i_v)_*i_v^*B$.
 By Lemma 4.10, the map from $B$ to $\prod_{v \in S} (i_v)_*i_v^*B$ induces an isomorphism on stalks, and hence is an isomorphism
by Lemma 4.9.

 Let $j = j_{K/F,S}$. By Theorem 4.8, Lemma 4.10b, and the fact that cohomology commutes
with direct products, we obtain

\proclaim{Corollary 4.12} 
 If $q > 0$, $H^p(\bar Y, R^qj_* \tilde A) = \coprod H^p(W_{\kappa(v)}, i_v^*R^qj_*\tilde A)$, where the
sum is taken over all $v \in S$. \endproclaim  

The next section will be devoted to computing these cohomology groups for small
values of
$p$ and $q$.

\proclaim {Proposition 4.13} The natural map $\phi$ from the sheaf $\Bbb Z$ on $T_{K,S}$ to the sheaf
$j_*j^*\Bbb Z$ is an isomorphism. \endproclaim

Proof. It is clear that $\phi$ is injective.  Let $\Cal X = (X_v)$ be an object of $T_{K,S}$. Let $f$ be in $j_*j^*\Bbb
Z (\Cal X) = j_*\tilde \Bbb Z(\Cal X) = Map (X_{v_0}, \Bbb Z) $  (see the remarks at the end of \S 1). Let $\{X_n\}$ be
the disjoint open cover of $X_{v_0}$ defined by $X_n = f^{-1}(n)$.  Let $X_{n,v} = f_v^{-1}(X_n)$.  Then the collection
$(X_{n,v})$ is a disjoint cover of $\Cal X$, and the element $g$ which takes each $(X_{n,v)}$ to $n$ lives in the \v Cech
cohomology of
$\Bbb Z$ with respect to this cover of
$\Cal X$, and its image in $\Bbb Z (\Cal X)$ maps to $f$.

\heading \S 5.  The computation of $H^p(\bar Y, R^qj_*\Bbb Z)$ and $H^p(\bar Y, R^qj_*\Bbb R)$ \endheading

\proclaim {Lemma 5.1}.  Let $G$ be a discrete group, and let $E$ be a sheaf on $T_G$.  Then the canonical map from $\widetilde
{E(G)}$ to $E$ induces an isomorphism of cohomology. \endproclaim.

Proof.  Since $G$ is discrete, any covering of a discrete $G$-space $X$ by $G$-spaces $X_i$ has a refinement consisting of
the $X_i$'s with the discrete topology.  It then follows by a standard comparison theorem in the theory of Grothendieck
topologies that the $T_G$-cohomology of any discrete $G$-space $X$ is the same as the cohomology of $X$ in the standard
topology of discrete $G$-sets and families of surjective morphisms.  But sheaves in this topology may be identified with
$G$-modules by making a sheaf $F$ correspond to the $G$-module $F(G)$.  ($G$ is a left $G$-space by left multiplication, and
the $G$-action on $F(G)$ is induced by letting $\sigma \in G$ act on $G$ by right multiplication by $\sigma^{-1}$).

Putting this together, we may identify the cohomology groups $H^p_{T_G}(*, E)$ with the groups $H^p(G, E(G))$, where the
cohomology is defined by the usual cochain definition.  

\proclaim{Lemma 5.2} Let $v$ be a finite place, let $j = j_{K,S}$, let $A$ be a continuous $W_F$-module, let $E =
i_v^*R^qj_* \tilde A$, and let
$G =
\Bbb Z = W_{\kappa(v)}$.  Then a)
$E(G) = H^q(\theta_v(I_v), A)$, and hence $H^p(W_{\kappa(v)}, i_v^*R^qj_* \tilde A) = H^p(W_{\kappa(v)}, H^q(\theta_v(I_v), A))$

b) $\Lim_{K,S} H^p(W_{\kappa(v)}, i_v^*R^qj_* \tilde A) = H^p(W_{\kappa(v)}, H^q(I_v, A))$.

\endproclaim

Proof.  Since $G$ has no non-trivial covers,
it is immediate that $i_v^pR^qj_*\Bbb Z(G)$ is naturally isomorphic to $E(G)$.   (Here $i_v^p$ denotes the presheaf inverse
image).  Recall that, if $C$ is a presheaf on $T_{K,S}$,  $i_v^p(C)(G)$ is given by the direct limit of those $C(\Cal X)$ for
which\ there is a map
$\phi: G
\to i_v^{-1}(\Cal X)$.  This then is equal to $H^q(T_{W_{K/F,S}}, W_{K/F,S    }/\theta_v(I_v), A)$. 

  By Theorem 4.6, this is just $H^q(\theta (I_v), A)$, and an application of Lemma 5.1 completes the proof of a).  Now
observe that $\Lim_{K.S} (H^q(W_{\kappa(v)}, H^q(\theta_v(I_v), A)) = $ (since
$W_{\kappa(v)} =
\Bbb Z)$
$H^p(W_{\kappa(v)}, \Lim_{K,S}H^q(\theta_v(I_v), A) = $(since $\theta(I_v)$ is compact and Moore cohomology commutes with limits
for compact groups) $H^p(\Lim_{K,S} \theta_v(I_v), A) = $ (by Lemma 3.1) $H^p(W_{\kappa(v)}, H^q(I_v, A))$, which shows b).

\proclaim {Lemma 5.3} Let $v$ be a finite place.  Then a) $H^1(I_v, \Bbb Z) =  H^1(I_v, \Bbb R) = 0$, 

 b) $H^0(W_{\kappa(v)},
H^2(I_v,
\Bbb Z))$ is naturally isomorphic to the Pontriagin dual $U_v^D$ of the local units $U_v$ in the completion $F_v$ of the field $F$
at
$v$, c) $H^0(W_{\kappa(v)}, H^2(I_v, \Bbb R)) = 0$.  
\endproclaim

Proof. If $A = \Bbb Z$ or $\Bbb R$,  $H^1(I_v, A) = Hom(I_v, A) = 0$.  From the exact sequence $ 1 \to I_v \to G_v \to
\hat
\Bbb Z
\to 1$, we get the Hochschild-Serre spectral sequence $H^p(\hat \Bbb Z, H^q(I_v, \Bbb Z)) \Rightarrow H^{p+q}(G_v, \Bbb Z)$. This
spectral sequence yields the short exact sequence $0 \to H^2(\hat \Bbb Z, \Bbb Z) \to H^2(G_v, \Bbb Z) \to H^0(\hat \Bbb Z,
H^2(I_v, \Bbb Z)) \to 0$.

By local class field theory $H^2(G_v, \Bbb Z)$ is naturally isomorphic to $Hom(F_v^*, \Bbb Q/\Bbb Z)$, so the above exact
sequence shows that $ H^0(\hat \Bbb Z, H^2(I_v, \Bbb Z))$ is naturally isomorphic to $Hom
(U_v, \Bbb Q/\Bbb Z)$ which (since $U_v$ is profinite) is the Pontriagin dual of $U_v$. But since $W_{\kappa(v)} = \Bbb Z$ is
dense in $\hat \Bbb Z$, 5.3b) follows immediately.  Since $I_v$ is compact, $H^2(I_v, \Bbb R) = 0$, which proves c).

\proclaim {Lemma 5.4} Let $\theta: H \to G$ be a map of topological groups, so we may regard any $G$-space as an $H$-space via
$\theta$. Let $I$ be a topological subgroup of $H$.   Let $Z$ be any topological space, regarded as a $G$-space with trivial
$G$-action, and let $X$ be any $G$-space.  Then any $H$-map $\phi: H/I \times Z$ to $X$ factors through the $G$-space
$G/\theta(I) \times Z$.  \endproclaim

Proof.  This follows immediately from the remark after Definition 4.2.

\proclaim {Lemma 5.5} Let $G$ be a connected topological group, and let $X$ be a topological space on which $G$
acts trivially.  Then $\check H^q(X, \Bbb Z)$ is naturally isomorphic to $\check H^q_{top}(X, \Bbb Z)$. \endproclaim

Proof. We first claim that any local-section $G$-cover $\rho_i: \{X_i \to X\}$ has a refinement by a cover of the form $\{G
\times V_i\}$, where
$\{V_i\}$ is an open cover of $X$, and $G$ acts on $G \times V_i$ by left multiplication on the first factor.  Given $x \in
X$, let $U_x$ be an open neighborhood of $x$ such that $s_x: U_x \to X_{i(x)}$ is a section of $\rho_{i(x)}$.  Define $\phi_x:
G \times U_x \to X_{i(x)}$ by $\phi_x(g, u) = gs_x(u)$, and verify first that $\phi_x$ is a $G$-map and next that $pr_2 =
\rho_i(x)\phi_x$, thus showing that $\{G \times U_x\}$ refines $\{X_i\}$.

We next claim that the \v Cech complex for the sheaf  
$\Bbb Z$  of the $G$-cover  $\{G \times V_i\}$ is the same as the \v
Cech complex of the cover $\{V_i\}$ of $X$.  This follows immediately because any map from a power $G^n$ of the connected group
$G$ to the discrete group $\Bbb Z$ is constant.

\proclaim {Lemma 5.6} Let $Z$ be a contractible topological space.  Let $v$ be a fixed archimedean place of $\bar Y_{K/F, S}$
and let $H = W_{\kappa(v)}$. Let $H$ act on $H \times Z$ by left multiplication on the first factor. We claim that 

a)$(i_v^pR^1j_*\Bbb Z) (H \times Z) = 0$. 

b) $(i_v^pR^2j_*\Bbb
Z)(H) = H^2(I_v, \Bbb Z)$. 

c) $(i_v^*R^2j_*\Bbb Z)(H) = H^2(I_v, \Bbb Z)$.

\endproclaim

Proof.  Let $E$ be any sheaf on $\bar Y_{K/F, S}$.  Then by definition, $(i_v^p(E))(H \times Z)$  is equal to the direct
limit of the  $E((X_w, f_w))$, where $H\times Z \to i_v^{-1}( (X_w, f_w)) = X_v$.  Now let $E = R^qj_*\Bbb Z$. It is
immediate that we may assume in the direct limit that
$X_v = H \times Z$ and that $X_w$ is the empty set if $w$ is neither $v$ nor the generic point $v_0$. Lemma 5.1 shows
that we may assume that $X_{v_0}$ is $W_{K/F}/\theta_v(I_v)$ and hence that $R^qj_*\Bbb Z((X_w, f_w)) =
H^q_{W_{K/F}}((W_{K/F}/(\theta_v(I_v)) \times Z), \Bbb Z)$.  By Lemma 4.3, this is the same as $H^q_{I_v}(Z, \Bbb Z)$. If
$q = 1$ this is equal to $\check H^1_{I_v}(Z, \Bbb Z)$ which in turn is equal by Lemma 5.5 to $\check H^1_{top}(Z, \Bbb Z)$
which is zero since $Z$ is contractible.

If $q =2$ and $Z$ is a point we have that $i_v^p(R^2j_*\Bbb Z) (H) = H^2_{I_v}(*, \Bbb Z) = H^2(I_v, \Bbb Z)$.  c) then
follows immediately because $H$ has no non-trivial covers.

\proclaim {Lemma 5.7} Let $G$ be a topological group, and $n$ a positive integer. Then $G^n$, regarded as a $G$-space with
$G$ acting diagonally, is isomorphic to $G \times G^{n-1}$ where $G$ acts by left multiplication on the first factor and
trivially on $G^{n-1}$. \endproclaim

Proof. Let $\phi: G^n \to G \times G^{n-1}$ by $\phi(g_1,\dots g_n) = (g_1, g_1^{-1}g_2, \dots g_{n-1}^{-1}g_n)$. It is easy
to see that $\phi$ is a $G$-isomorphism.

\proclaim {Proposition 5.8} Let $v$ be an archimedean place.  

a) $H^p(W_{\kappa(v)}, i_v^*R^1j_* \Bbb Z) = 0$ for $p = 0, 1, $ and $2$.

a') $H^p(W_{\kappa(v)}, i_v^*R^1j_* \Bbb R) = 0$ for $p = 0, 1, $ and $2$.

b) $H^0(W_{\kappa(v)}, i_v^*R^2j_*\Bbb Z) = H^2(I_v, \Bbb Z)^{W_{\kappa(v)}}$	 

b') $H^0(W_{\kappa (v)}, i_v^*R^2j_*\Bbb R) = 0$.

c) $H^2(I_v, \Bbb Z)^{W_{\kappa(v)}} = {U_v}^D$.

\endproclaim

Proof. We have the standard spectral sequence from \v Cech to derived functor cohomology:

$$ E^{p.q}_2 = \check H^p(\kappa (v), \underline{H}^q(i_v^*R^1j_*\tilde A)) \Rightarrow H^{p+q}(\kappa(v), i_v^*R^1j_*\tilde A)$$
where we know that $E^{0,q}_2 = 0$ for $q > 0$.

We begin with the case $p = 2$.  The spectral sequence immediately gives the exact sequence;

$$0 \to \check H^2(\kappa(v), i_v^*R^1j_*\tilde A) \to H^2(\kappa(v), i_v^*R^1j_*\tilde A) \to \check H^1(\kappa(v), \underline
{H}^1(i_v^*R^1j_*\tilde A)) $$

So it suffices to show that the first and third terms in this exact sequence are zero.  We begin with the first: 

We first let $A = \Bbb Z$ and show that, more generally, $\check H^p(\kappa(v), i_v^*R^1j_*\Bbb Z) = 0$.  Since the covering
$\{H\}$ of
$*$ is initial, it is enough to show that $(i_v^*R^1j_*\Bbb Z )(H^n) = 0$.  By lemma 5.7, this is equivalent to showing that
$(i_v^*R^1j_*\Bbb Z) (H \times H^{n-1}) = 0$, where $H$ acts trivially on $H^{n-1}$. But this is an immediate consequence of Lemma
1.6 and Lemma 5.5, since $H$ is contractible and locally contractible.

Now let $A = \Bbb R$.  If $E$ is any sheaf of $\Bbb R$ -vector spaces on $H \times H^{n-1}$, and $q > 0$, Corollary 1.7 shows that
$H^q _H(H\times H^{n-1}, E)$ is isomorphic to  $H^q(H^{n-1}, \alpha_*E)$  which is equal to zero by Lemma 2.5.

Now we look at the third term.  Since $H^1 = \check H^1$, we have to show that $\check H^1(H \times H^{n-1}, i_v^*R^1j_*\tilde
A) = 0$.  A typical term in a coinitial cover of $H \times H^{n-1}$ is $H^r \times X$ with $X$ contractible, locally
contractible, and metrizable. But rewriting this as $H \times (H^{r-1} \times X)$ and again using Lemma 5.5 in the case when $A =
\Bbb Z$ and Lemma 2.5 when $A = \Bbb R$ enables us to copy the arguments of the preceding paragraph, since $H^{r-1} \times X$ is
also contractible, locally contractible, and metrizable.

The case when $p =1$ is similar but easier.

Now b) and b') follow immediately from Lemma 5.6c.

If $v$ is complex $H^2(I_v, \Bbb Z) = \Bbb Z$, and if v is real we have the exact sequence $1 \to S^1 \to I_v \to \Bbb Z/2\Bbb Z
\to 1 $.  The Hochschild-Serre spectral sequence shows that $H^2(I_v, \Bbb Z) = H^2(S^1, \Bbb Z)^{\Bbb Z/2\Bbb Z} = \Bbb Z/2\Bbb
Z$.  In both cases the Weil group $W_{\kappa(v)}$ acts trivially, and we get $\Bbb Z$ and $\Bbb Z/2\Bbb Z$, the duals of $S^1$
and $\pm 1$ respectively.

\proclaim{Theorem 5.9} Let $A$ be either $\Bbb Z$ or $\Bbb R$. a) $H^p(\bar Y_{K,S}, R^1(j_{K,S})_* \tilde A) =0$ for $p = 0, 1,
2$.

b) $H^p(\bar Y, R^1j_* \tilde A) = 0$ for $p = 0, 1, 2$.

c) $H^0(\bar Y_{K,S} , R^2(j_{K,S})_* \Bbb Z) = \coprod _{v \in S} (U_v)^D$. 

c') $H^0(\bar Y, R^2(j_{K,S})_* \Bbb R) = 0$.

d) $H^0(\bar Y, R^2j_* \Bbb Z) = \coprod_{v \neq v_0} (U_v)^D$.

d') $H^0(\bar Y, R^2j_*\Bbb R) = 0$.
 \endproclaim

Proof.  Parts a) and c) follow immediately from Corollary 4.12, Lemma 5.2, Lema 5.3, and Proposition 5.8.  Parts b) and d) follow
from a) and c) by taking limits.

Let $Pic(\bar Y)$ be the Arakelov class group of $F$, i. e., the group obtained by taking the id\`ele group of $F$ and dividing
by the principal id\`eles and the unit id\`eles (a unit id\`ele $(u_v)$ is defined by $|u_v|_v = 1$ for all $v$).  Let $Pic^1(\bar
Y)$ be the kernel of the absolute value map from $Pic (\bar Y)$ to $\Bbb R^*$. Let $\mu(F)$ denote the group of roots of unity in
$F$.

\proclaim{Theorem 5.10} a) $H^0(\bar Y, \Bbb Z) = \Bbb Z$ 

b) $H^1(\bar Y, \Bbb Z) = 0$

c) $H^2(\bar Y, \Bbb Z) = (Pic^1(\bar Y)^D$

d) $H^3(\bar Y, \Bbb Z) = \mu(F)^D$. \endproclaim

Proof. a) is clear.  The Leray spectral sequence for $j_*$ gives first that $H^1(\bar Y, \Bbb Z) = H^1(W_F, \Bbb Z) = 0$, which
proves b). Next it gives (using Theorem 5.9) the exact sequence

$$ 0 \to H^2(\bar Y, \Bbb Z) \to H^2(W_F, \Bbb Z) \to \coprod_{v \neq v_0} (U_v)^D \to H^3(\bar Y, \Bbb Z) \to H^3(W_F, \Bbb Z) =
0 $$.

This is easily seen (using Theorem 3.6) to be the Pontriagin dual of the sequence:

$$ 0 \to H^3(\bar Y, \Bbb Z)^D \to \prod _{v \neq v_0} U_v \to C^1_F \to H^2(\bar Y, \Bbb Z)^D \to 0$$

which completes the proof, since the roots of unity are the  kernel of the map from the unit id\`eles to $C^1(F)$ and $Pic^1(F)$ is
defined to be the cokernel.

\proclaim {Theorem 5.11} a) $H^0(\bar Y, \Bbb R) = \Bbb R$

b) $H^1(\bar Y, \Bbb R) = H^2(\bar Y, \Bbb R) = 0$.  \endproclaim

Proof. a) is clear, and b) follows from the Leray spectral sequence, using Theorem 5.9.

\heading \S 6. Cohomology with compact support \endheading

Let $Y$ be Spec $O_F$ abd let $\varphi$ be the natural inclusion of $Y$ in $\bar Y$.  Let $E$ be any sheaf on $Y$.  We define the
sheaf
$\varphi_!E$ on
$\bar Y$ to be the sheaf associated with the presheaf $P$ defined by $P(\Cal X = (X_v)) = E((X_v))$ if $X_v = \phi$ for all $v$ not
in $Y$, and
$P(\Cal X) = 0$ otherwise. 
\proclaim {Proposition 6.1} Let $F$ be any sheaf on $\bar Y$. There exists an exact sequence of sheaves on $\bar Y$:

$$ 0 \to \varphi_! \varphi^*F  \to F \to i_*i^*F \to 0 $$

where $i_*i^*F = \prod_{v \in Y_{\infty}}(i_v)_*i_v^*F$.

\endproclaim

Proof. We first show that for all $v$ in $\bar Y$, there exists an exact  sequence 

$$ 0 \to i_v^* \varphi_! \varphi^*F \to i_v^* F \to i_v^*i_* i^*F \to 0 $$

We first see easily that if $v$ is non-archimedean that $i_v^*\varphi_! \varphi^*F = i_v^*F$, and $i_v^*(i_*i^*F) = 0$ by Lemma
4.10c), so we get exactness.  If $v$ is archimedean, $i_v^*\varphi _! \varphi^*F = 0$ and $i_v^* (i_*i^*F) = i_v^*i_*F$ by
Lemma 4.10b), so again we get exactness.

The exactness of the above exact sequences implies the Proposition, using Lemma 4.9 and the fact that $i_v^*$ is exact
(Lemma 4.10 a)).  

\proclaim {Lemma 6.2} Let $v$ be an archimedean valuation.  Then a) $H^i(W_{\kappa(v)}, \Bbb Z) = 0$ for $i > 0$. 

b) $H^i(\bar Y, i_*\Bbb Z) = 0$ for $i > 0$.

\endproclaim

Proof. a) is immediate because $W_{\kappa(v)} = \Bbb R$, $\Bbb R$ is contractible, and $\Bbb Z$ is discrete.  Then b) follows
becaue $i_*$ is exact.

\proclaim {Theorem 6.3}a) $H^0(\bar Y, \varphi_! \Bbb Z) = 0$

b) $ H^!(\bar Y, \varphi_! \Bbb Z) = (\coprod_{S_\infty} \Bbb Z)/\Bbb Z)$

c) $H^2(\bar Y, \varphi_! \Bbb Z) = Pic^1(\bar Y)^D$

d) $H^3(\bar Y, \varphi_! \Bbb Z) = \mu(F)^D$ \endproclaim

Proof.  This follows immediately from Theorem 5.10, Proposition 6.1, and Lemma 6.2.

\proclaim{Proposition 6.4} There is a natural exact sequence

$$ 0 \to Pic(Y)^D \to Pic^1(\bar Y) ^D \to Hom(U_F, \Bbb Z) \to 0 $$ \endproclaim

Proof.  Let $F_v$ denote the completion of $F$ at the archimedean valuation $V$. Then we have a natural inclusion $i$ of $\prod_v
F_v^*$ into the id\`ele group $J_F$.  We then obtain an exact sequence;

$$ \CD 0 @ >>> \coprod \Bbb R^*_{>0} @ >\tilde i>> Pic(\bar Y) @ >>> Pic(Y) @>>> 0 \endCD $$

where $\tilde i$ is induced by $i$.

Then the logarithmic embedding of the units yields the exact sequence

$$ 0 \to V/L \to Pic^1(\bar Y) \to Pic(Y) \to 0 $$

where $V$ is the kernel of the sum map from $\coprod_v \Bbb R$ to $\Bbb R$, and $L$ is the lattice in $V$ obtained by
taking the image of the unit group $U_F$ under the map which sends a unit $u$ to the vector $(log|u|_v)$.

We now examine the following commutative diagram:

$$ \CD
0 @ >>> (V/L)^D @>>> V^D @ >>> L^D @ >>> 0 \\
@ . @ AAA @ A\alpha AA @ AAA \\
@ . 0 @ >>> Hom (V, \Bbb R) @>\beta >> Hom (L. \Bbb R) \\
@ .@ . @ AAA @ A\gamma AA\\
@ .@ . 0 @ >>> Hom (L, \Bbb Z) \\
\endCD $$
where $\alpha$ and $\beta$ are isomorphisms, and $\gamma$ is injective.  This defines an isomorphism between $Hom(L, \Bbb Z)$ and
$(V/L)^D$, and the proposition follows, after we observe that the natural map from $Hom(L, \Bbb Z)$ to $Hom (U_F, \Bbb Z)$ is an
isomorphism.

\heading \S 7. Euler characteristics \endheading

Let $n \geq 1$ and let  $$ \CD 
0 @ >>> V_0 @ >T_0>> V_1 @>T_1>> \dots @>T_{n-1} >> V_n @ >>> 0 \\
\endCD $$

be an exact sequence of real vector spaces, and let $B_i$ denote an ordered basis for $V_i$. We recall the definition of the
determinant of the above data.  If $n = 1$ the data determine an $n\times n$ matrix, and we take the determinant of that matrix.

If $n = 2$, let $B_0 = (a_1, \dots a_r)$, let $B_1 = (b_1, \dots b_{r+s})$, and let $B_2 = (c_{r+1}, \dots c_{r+s})$. For $1\leq
i \leq r$, let $d_i = T_0(a_i)$.  Choose $(d_{r+1}, \dots d_{r+s})$ in $V_1$ such that $T_1(d_i) = c_i$.  In the one-dimensional
space $\Lambda^{r+s}V_1$ the element $d_1 \wedge d_2 \wedge \dots d_{r+s}$ is clearly independent of the choice of $d_i$, and
we define our determinant $\delta$ so that $d_1 \wedge d_2 \dots \wedge d_{r+s}$ = $\delta (b_1 \wedge b_2 \wedge \dots
\wedge b_{r+s})$.

We finish by giving an inductive definition. Assume we have defined the determinant for $n \leq N$  and we wish to define it for
$n = N+1$. We let $I$ be the image of $T_{n-1}$ so that we have the two exact sequences :

$$ \CD
0 @ >>> V_0 @ >T_0>> V_1 @ >T_1>> \dots @ >T_{n-1} >> I @ >>> 0 \\
\endCD $$

$$ \CD
0 @ >>> I @ >i>> V_n @ >T_n>> V_{n+1} @ >>> 0 \\
\endCD $$

where $i$ is the inclusion of $I$ in $V_n$.  We choose any basis $C$ for $I$. We now define the determinant $\delta$ of the
sequence
$$ \CD
0  @ >>> V_0 @ >T_0>> \dots @ >T_n>> V_{n+1} @ >>> 0 \\
\endCD  $$

with bases $B_0, \dots, B_{n+1}$, to be $\delta_1 (\delta_2)^{(-1)^n}$,
where $\delta_1$ is the determinant of the sequence 

$$ \CD
0 @ >>> V_0 @ >T_0>> V_1 @ >T_1>> \dots @ >T_{n-1}>> I @ >>> 0 \\
\endCD $$ 

where $V_i$ has basis $B_i$ and $I$ has basis $C$,
and $\delta_2$ is the determinant of the sequence

$$ \CD
0 @ >>> I @ >i>> V_n @ >T_n>> V_{n+1} @ >>> 0 \\
\endCD $$

where $I$ has basis $C$, and $V_n$ and $V_{n+1}$ have bases $B_n$ and $B_{n+1}$. 
It is easy to see that this definition is independent of the choice of $C$.

Now let $A_0$, $A_1$, $\dots$, $A_n$ be finitely generated abelian groups, and let $V_i = A_i \otimes \Bbb R$. Assume that there
exist $\Bbb R$-linear transformations $T_i:V_i \to V_{i+1}$ such that the sequence

$$ \CD
0 @ >>> V_0 @ >T_0>> V_1 @ >T_1>> \dots  V_{n-1} @ >T_n>> V_n @ >>> 0 \\
\endCD $$

is exact.

We define the Euler characteristic $\chi (A_0, A_1,\dots A_n, T_0 \dots T_{n-1})$ to be the alternating product $\prod_{i =0}^n
|((A_i)_{tor})|^{(-1)^i}$ divided by the determinant of $(V_0, \dots V_n,  T_0, \dots T_{n-1}, B_0, \dots, B_n)$, where the
$B_i$ are the images of bases of the free abelian groups $A_i/(A_i)_{tor}$.  

The $B_i$ are of course not unique, but a change of basis only changes the determinant by the determinant of a matrix in $GL(r,
\Bbb Z)$, i. e. by $\pm 1$.

So our Euler characteristic is well-defined up to sign.

\heading \S 8. Dedekind zeta-functions at zero \endheading

In this section we wish to verify that the conjecture stated in the introduction is true for Dedekind zeta-functions, modulo the
assumption that the higher cohomology groups are zero.

We first define our Euler characteristic.    Let $F$ be a number field, let $O_F$ be the ring of integers in $F$, and let $Y =$
Spec $O_F$.  Let $\bar Y$ be $Y$ together with the archimedean primes of $F$, given the Weil -\'etale topology as above.  Let
$\varphi$ be the inclusion of $Y$ in $\bar Y$.   

Let $\Bbb R$ denote the sheat on $\bar Y$ determined by defining $\Bbb R((X_v))$ to be the set of compatible continuous $W_v$-maps
from $X_v$ to $\Bbb R$, where $W_v$ acts trivially on $\Bbb R$.  It is clear both that this is a sheaf and that such a set is
determined  by giving a $W_{v_0}$-map from $X_{v_0}$ to $\Bbb R$. It is also clear that this is the same sheaf as the sheaf
$j_*\tilde \Bbb R$.

The Leray spectral sequence for the map $j_*$ yields:
$$ 0 \to H^1(\bar Y, \Bbb R) \to H^1(W_F, \Bbb R) \to H^0(\bar Y, R^1j_*\Bbb R) \to H^2(\bar Y, \Bbb R) \to H^2(W_F, \Bbb R) = 0 $$

where $H^2(W_F, \Bbb R) = 0$ by Lemma 3.4.

But $R^1j_*\Bbb R$ is isomorphic to $\coprod (i_v)_*i_v^*R^1j_*\Bbb R$, and so we conclude easily that $H^1(\bar Y, R^1j_*\Bbb R)$
is isomorphic to $\coprod H^1(I_v, \Bbb R)$, where the sums are taken over all non-trivial  valustions of $F$.  But whether $v$
is archimedean or non-archimedean, $I_v$ is compact, so $H^1(I_v, \Bbb R) = Hom (I_v, \Bbb R) = 0$. We conclude that $H^1(\bar Y,
\Bbb R) = H^1(W_F, \Bbb R) = Hom (W_F, \Bbb R)$, and that $H^2(\bar Y, \Bbb R) = 0$.  Let $\psi$ in $H^1(\bar Y, \Bbb R)$ be the
homomorphism obtained by mapping
$W_F$ to its abelianization $C_F$ and then taking the logarithm of the absolute value.  

We next observe that first, by standard arguments the category of sheaves of $\Bbb R$-modules has enough injectives, and next,
that any injective sheaf of $\Bbb R$-modules is injective as a sheaf of abelian groups.  These observations imply that taking the
Yoneda product with $\psi$ in $H^1(\bar Y, \Bbb R) = Ext^1_{\bar Y}(\Bbb R, \Bbb R$) induces a map from $H^q(\bar Y, F) =
Ext^q_{\bar Y}(\Bbb R, F)$ to
$H^{q+1}(\bar Y, F) = Ext^{q+1}_{\bar Y}(\Bbb R, F)$, where $F$ is any sheaf of $R$-modules.

\proclaim{Theorem 8.1} Assume that $H^q(\bar Y, \varphi_! \Bbb Z) = 0$ for $q > 3$.  Let $\zeta_F$ be the Dedekind
zeta-function of $F$.  Then the Euler characteristic  $\chi (H^*(\bar Y, \varphi_! \Bbb Z))$ is well-defined and is equal to
$\pm \zeta_F^*(0)$.
\endproclaim

Proof.  We first observe that the groups $H^i(\bar Y, \varphi_!\Bbb Z)$ are finitely-generated by Theorem 6.3 and Proposition 6.4. 
We must show next that the natural map from
$H^i(\bar Y, \varphi_! \Bbb Z) \otimes \Bbb R$ to $H^i(\bar Y, \varphi_! \Bbb R)$ is an isomorphism. Look at the commutative
diagram:
$$ \CD
@ . H^2(\bar Y, \Bbb R) = 0 \\
@. @ AAA \\
@  . H^2(\bar Y, \varphi_!\Bbb R) \\
 @. @  AAA  \\
0 @>>> H^1(\bar Y, i_*\Bbb R)@ >\alpha  >> H^1(\bar Y, i_*S^1)@ >>> H^2(\bar Y, i_* \Bbb Z) = 0 \\
@ . @ A\gamma AA @ AAA @ AAA \\
@  . H^1(\bar Y, \Bbb R) @ >>> H^1(\bar Y, S^1) @ >>> H^2(\bar Y, \Bbb Z) @>>> H^2(\bar Y, \Bbb R) = 0 \\
 @ . @ A\delta AA @ AAA @ A\beta AA @ AAA \\
 @ . H^1(\bar Y, \varphi_! \Bbb R)  @ >>> H^1(\bar Y,\varphi_!S^1) @ >>> H^2(\bar Y, \phi_!\Bbb Z) @ >\epsilon >> H^2(\bar Y,
\varphi_!\Bbb R)
\\
\endCD $$

It is easy to see that $\gamma$ is injective, so $\delta$ is the zero map, so $H^1(\bar Y, \varphi_!\Bbb R)$ may be identified with
$H^1(\bar Y, \varphi_!\Bbb Z) \otimes \Bbb R$ , and we take a basis of $H^1(\bar Y, \varphi_! \Bbb R)$ obtained by choosing $r_1 +
r_2 -1$ archimedean primes of $F$.

By a tedious but straightforward calculation with injective resolutions, we see that the map $\epsilon$ may be computed by
applying $\beta$, lifting to $H^1(\bar Y, S^1)$, mapping to $H^1(\bar Y, i_*S^1)$, applying $\alpha ^{-1}$, and mapping to
$H^2(\bar Y, \varphi_! \Bbb Z)$. 

Now by comparing this diagram with the diagram at the end of Section 6, we see that we may first identify $H^2(\bar Y, \varphi_!
\Bbb R)$ with $Hom(V_0, \Bbb R)$, where $V = \coprod_{v \in S_\infty} \Bbb R^*_{>0}$ and $V_0$ is the kernel of the product map to
$\Bbb R^*_{>0}$.  Next, we may take as a basis of this group coming from $H^2(\bar Y, \varphi_! \Bbb Z)$ the dual basis of any
basis for the units of $F$ modulo torsion, identifying $V_0$ with $U_F \otimes \Bbb R^*_{>0}$ via the map $u \mapsto ( |u|_v)$ for
the same set of
$r_1 +r_2 -1$ $v's$ we used above.
Finally the Yoneda product with $\psi$ clearly takes $1_v$ to the map $f_v$ where $f_v ((x_w)) = \log x_v$.

It is now easy to see that the determinant of the pair consisiting of $H^*(\bar Y, \varphi_! \Bbb Z)$ and Yoneda product with
$\psi$ is
$R^{-1}$ where $R$ is the classical regulator.

Since $H^0(\bar Y, \varphi_! \Bbb Z) = 0$, $(H^1(\bar Y, \varphi_! \Bbb Z))_{tor} = 0$. $|(H^2(\bar Y, \varphi_! \Bbb Z))_{tor}| =
h$, and $|H^3(\bar Y, varphi_! \Bbb Z)| = w$, the Euler characteristic of $H^*(\bar Y, \varphi_!, \Bbb Z)$ is equal to $hR/w$ which
up to sign is $\zeta^*_F(0)$.

\heading References \endheading

\noindent[A-T] Artin, E. and Tate, J.  Class Field Theory, Benjamin, New York, 1967

\noindent [A] Artin, M. Grothendieck Topologies, mimeographed notes, Harvard University, 1962

\noindent [B-W] Borel, A., and Wallach, N. Continuous cohomology discrete subgroups, and representations of reductive groups,
Annals of Mathematics Studies, Princeton University Press, Princeton, 1980 

\noindent [Ge1] Geisser, T. Weil-\'etale cohomology over finite fields (to appear, Math. Annalen)

\noindent [Ge2] Geisser, T. Arithmetic cohomology over finite fields and special values of $\zeta$-functions, preprint 2004

\noindent [Go] Godement, R. Topologie alg\'ebrique et th\'eorie des faisceaux, Hermann, Paris, 1958 

\noindent [L] Lichtenbaum, S. The Weil-\'etale topology on schemes over finite fields. (to appear, Compositio Math.)

\noindent [M1] Moore, C. C. Extensions and low dimensional cohomology of locally compact groups I, II, Trans. Amer. Math. Soc. 113
(1964) 40-63, 64-86

\noindent [M2] Moore, C. C.  Group extensions and cohomology for locally compact groups III, Trans.  Amer. Math. Soc. 221 (1976)
1-33

\noindent [NSW] Neukirch, J., Schmidt, A., and Wingberg, K., Cohomology of Number Fields, Springer, New York, 2000

\noindent [R] Rajan, C. S. On the Vanishing of the Measurable Schur Cohomology Groups of Weil Groups, Compositio Math. 140 (2004)
84-98

\noindent [T] Tate, J., Number theoretic background, Proc. Sump. Pure Math., vol. 33, Part 2, Amer. Math. Soc. Providence, R. I.
1977, 3-26

\noindent [W] Wigner, D.  Algebraic cohomology of topological groups, Trans. Amer. Math. Soc. 178 (1973) 83-93

\bye